\documentclass[graybox,footinfo]{svmult}

\smartqed
\usepackage{mathptmx}       
\usepackage{helvet}         
\usepackage{courier}        
\usepackage{type1cm}        
\usepackage{graphicx}       

\usepackage{array,colortbl}
\usepackage{amsmath,amsfonts,amssymb,bm} 
\DeclareFontFamily{U}{mathx}{\hyphenchar\font45}
\DeclareFontShape{U}{mathx}{m}{n}{
      <5> <6> <7> <8> <9> <10>
      <10.95> <12> <14.4> <17.28> <20.74> <24.88>
      mathx10
      }{}
\DeclareSymbolFont{mathx}{U}{mathx}{m}{n}
\DeclareFontSubstitution{U}{mathx}{m}{n}
\DeclareMathAccent{\widecheck}      {0}{mathx}{"71}

\renewcommand{\email}[1]{\emailname: #1} 

\usepackage{mathptmx}       
\usepackage{helvet}         
\usepackage{courier}        

\usepackage{makeidx}         
\usepackage{graphicx}        
\usepackage[bottom]{footmisc}

\usepackage{latexsym}
\usepackage{amsmath}
\usepackage{amsfonts}
\usepackage{amssymb}
\usepackage{bm}

\usepackage{url}
\usepackage{algorithm}
\usepackage{algorithmic}
\usepackage[misc,geometry]{ifsym}

\spdefaulttheorem{assumption}{Assumption}{\upshape \bfseries}{\itshape}
\spdefaulttheorem{algo}{Algorithm}{\upshape \bfseries}{\itshape}

\DeclareSymbolFont{bbold}{U}{bbold}{m}{n}
\DeclareSymbolFontAlphabet{\mathbbold}{bbold}

\DeclareSymbolFont{bbold}{U}{bbold}{m}{n}
\DeclareSymbolFontAlphabet{\mathbbold}{bbold}

\usepackage{microtype} 

\usepackage[colorlinks=true,linkcolor=black,citecolor=black,urlcolor=black]{hyperref}
\usepackage{bookmark}
\makeatletter 
  \providecommand*{\toclevel@author}{999}
  \providecommand*{\toclevel@title}{0}
\makeatother

\begin{document}

\title*{Estimates For Logarithmic and Riesz Energies For Spherical $t$-designs}
\author{Tetiana A. Stepanyuk}

\institute{Tetiana A. Stepanyuk 
\at $^{(1)}$Graz University of Technology, Kopernikusgasse 24, Graz, Austria;
$^{(2)}$Institute of Mathematics of Ukrainian National Academy of Sciences, 3, Tereshchenkivska st., 01601, Kyiv-4, Ukraine
\email{tania$_{-}$stepaniuk@ukr.net},
%
%
}
\maketitle


\abstract{In this paper we find asymptotic equalities for the  discrete logarithmic energy of sequences of well separated spherical $t$-designs on the unit sphere ${\mathbb{S}^{d}\subset\mathbb{R}^{d+1}}$, $d\geq2$. 
Also we establish exact order estimates for  discrete Riesz $s$-energy, $s\geq d$, of sequences of well separated spherical $t$-designs.
}

\section{Introduction}\label{intro} 

Let $\mathbb{S}^{d}=\{\mathbf{x}\in\mathbb{R}^{d+1}: \ |\mathbf{x}|=1\}$, where $d\geq2$, be the unit
sphere in the Euclidean space $\mathbb{R}^{d+1}$, equipped with the Lebesgue measure 
$\sigma_{d}$ normalized by $\sigma_{d}(\mathbb{S}^{d})=1$. 

\begin{definition}\label{def1}
  A spherical $t$-design is a finite subset $X_{N}\subset \mathbb{S}^{d}$ with
  a characterising property that an equal weight integration rule with nodes
  from $X_{N}$ integrates all spherical polynomials $p$ of total degree at most  $t$ exactly;
  that is,
\begin{equation*}
  \frac{1}{N}\sum\limits_{\mathbf{x}\in X_{N}}p(\mathbf{x})=
  \int_{\mathbb{S}^{d}}p(\mathbf{x})d\sigma_{d}(\mathbf{x}), \quad
  \mathrm{deg}(p)\leq t.
\end{equation*}

Here $N$ is the cardinality of $X_{N}$ or the number of points of spherical design.
\end{definition}

The concept of spherical $t$-design was introduced  by Delsarte, Goethals and
Seidel in the groundbreaking paper \cite{Delsarte-Goethals-Seidel1977:spherical_designs}, since then they attracted a lot of interest from scientific community (see   e.g.,   \cite{BrauchartGrabner}).

The logarithmic energy of a set of $N$ distinct points (or an $N$-point set) $X_{N}$ on 
$\mathbb{S}^{d}$ is defined as
\begin{equation}\label{RieszbDef}
  E_{log}^{(d)}(X_{N}):=
  {\mathop{\sum}\limits_{i,j=1,\atop i\neq j}^{N}}
  \log \frac{1}{|\mathbf{x}_{i}-\mathbf{x}_{j}|}=2 \sum\limits_{1\leq j<i\leq N}
  \log \frac{1}{|\mathbf{x}_{i}-\mathbf{x}_{j}|}.
\end{equation}

This paper investigates the logarithmic energy for spherical $t$-designs.
Spherical $t$-designs of a fixed strength $t$ can have points arbitrary close together (see, e.g. \cite{HesseLeopardiTheCoulombEnergy}), hence the logarithmic energy of $N$-point spherical $t$-designs  can have no asymptotic bounds in terms of $t$ and $N$.  That's why we will have additional condition and  consider the sequences of well--separated spherical $t$-designs.
\begin{definition}
  A sequence of $N$-point sets $X_{N}$,
  $X_{N}=\big\{\mathbf{x}_{1},\ldots, \mathbf{x}_{N} \big\}$, is called
  well-separated if there exists a positive constant $c_{1}$ such that
\begin{equation}\label{wellSeparat}
  \min\limits_{i\neq j}|\mathbf{x}_{i}- \mathbf{x}_{j}|>
  \frac{c_{1}}{N^{\frac{1}{d}}}.
\end{equation}
\end{definition}

The existence of  $N$-point spherical $t$-designs with $ N(t)\asymp t^{d}$ was proven by Bondarenko, Radchenko and Viazovska
\cite{Bondarenko-Radchenko-Viazovska2013:optimal_designs}.
They showed that for ${d\geq 2}$, there exists a constant $c_{d}$, which depends only of $d$, such that for every $N\geq c_{d}t^{d}$ there exists a spherical $t$-design on $\mathbb{S}^{d}$ with $N$ points. Two years later by these authors 
 in 
\cite{Bondarenko-Radchenko-Viazovska2015:Well_separated} the existence of  $N$-point well--separated spherical $t$-designs with $ N(t)\asymp t^{d}$ was proven.
Namely,
they showed that for each ${d\geq 2}$, $t\in \mathbb{N}$, there exist positive constants  $c_{d}$ and $\lambda_{d}$, depending only on $d$, such that for every $N\geq c_{d}t^{d}$, there exists a spherical $t$-design on $\mathbb{S}^{d}$,
consisting of $N$ points $\{\mathbf{x}_{i}\}_{i=1}^{N}$ with $|\mathbf{x}_{i}-\mathbf{x}_{j}|\geq \lambda_{d}N^{-\frac{1}{d}}$ for $i\neq j$.

On the basis of these results we always assume that  $ N=N(t)\asymp t^{d}$.

We write
 $a_{n}\asymp b_{n}$ to mean that there exist positive constants $C_{1}$ and
 $C_{2}$ independent of $n$ such that $C_{1}a_{n}\leq b_{n}\leq C_{2}a_{n}$ for
 all $n$.

 Denote by $\mathcal{E}_{log}^{(d)}(N)$ the minimal discrete logarithmic energy for $N$-points on the sphere
\begin{equation}\label{minLogEnergie}
 \mathcal{E}_{log}^{(d)}(N):=\inf\limits_{X_{N}}E_{log}^{(d)}(X_{N}),
\end{equation}
where the infimum is taken over all $N$-points subsets of $\mathbb{S}^{d}$.

From the papers of Wagner \cite{Wagner}, Kuijlaars and Saff \cite{KuijlaarsSaff:1998Asymptotics}  and  Brauchart \cite{Brauchart2008} it follows that for  $d\geq2$ and  as $N\rightarrow\infty$ the following asymptotic equality holds
 \begin{equation}\label{theoremLog}
\mathcal{ E}_{log}^{(d)}(N)=N^{2}\int\limits_{\mathbb{S}^{d}}\int\limits_{\mathbb{S}^{d}}\log\frac{1}{|\mathbf{x}-\mathbf{y}|}d\sigma_{d}(\mathbf{x})d\sigma_{d}(\mathbf{y})-\frac{1}{d}N\log N +\mathcal{O}(N).
\end{equation}
 
 Also in \cite{BoyvalenkovDragnevHardinSaffStoyanova} some general upper and lower bounds for the energy of spherical designs were found.

 We show that for every well-separated sequence of
 $N$-point spherical $t$-designs on 
   $\mathbb{S}^{d}$, $d\geq2$,   with $N\asymp t^{d}$ the following asymptotic equality holds
 \begin{equation*}
 E_{log}^{(d)}(X_{N})=N^{2}\int\limits_{\mathbb{S}^{d}}\int\limits_{\mathbb{S}^{d}}\log\frac{1}{|\mathbf{x}-\mathbf{y}|}d\sigma_{d}(\mathbf{x})d\sigma_{d}(\mathbf{y})  -\frac{1}{d}N\log N +\mathcal{O}(N).
\end{equation*}
 
Comparing two last formulas, we have that the leading and second terms are exactly the same, and third terms are of the same order. So, we can summarize, that for logarithmic energy well-separated spherical $t$-designs are as good  as point sets which  minimize the logarithmic energy.

For given $s>0$ the discrete Riesz $s$-energy of a set of $N$ distinct points (or an $N$-point set) $X_{N}$ on 
$\mathbb{S}^{d}$ is defined as
\begin{equation}\label{RieszbDef}
  E_{s}^{(d)}(X_{N}):=\frac{1}{2}
  {\mathop{\sum}\limits_{i,j=1,\atop i\neq j}^{N}}|\mathbf{x}_{i}-\mathbf{x}_{j}|^{-s},
\end{equation}
where $|\mathbf{x}|$ denotes the Euclidian norm in $\mathbb{R}^{d+1}$ of the vector $\mathbf{x}$. In the case $s=d-1$ the energy (\ref{RieszbDef}) is called as Coulomb energy.

Hesse \cite{Hesse:2009s-energy}
 showed,  that if spherical $t$-designs with 
${N=\mathcal{O}(t^{2})}$ exist, then they have asymptotically minimal Riesz energy $E_{s}^{(2)}(X_{N})$ for $s\geq 2$.
In particular, under the assumption that $N\leq \kappa t^{2}$, it was shown that for $s> 2$, there exists a positive constant $c_{s}$ such that for every  well separated sequence $N$
-point spherical $t$-designs the following estimate holds
\begin{equation}\label{Hesse}
  E_{s}^{(2)}(X_{N})\leq c_{s}N^{1+\frac{s}{2}},
\end{equation}
and for $s=2$, there exists a positive constant $c_{2}$, such that
\begin{equation}\label{Hesse_s=2}
  E_{s}^{(2)}(X_{N})\leq \frac{\sum\limits_{k=0}^{t}\frac{1}{k+1}}{4}N^{2}+c_{2}N^{2},
\end{equation}
and
\begin{equation}\label{Hesse_limit}
 \lim\limits_{N\rightarrow\infty}\frac{ E_{s}^{(2)}(X_{N})}{N^{2}\log N}=\frac{1}{8}.
\end{equation}

Denote by $\mathcal{E}_{s}^{(d)}(X_{N})$ the minimal discrete $s$-energy for $N$-points on the sphere
\begin{equation}\label{minCoulomb}
 \mathcal{E}_{s}^{(d)}(N):=\inf\limits_{X_{N}}E_{s}^{(d)}(X_{N}),
\end{equation}
where the infimum is taken over all $N$-points subsets of $\mathbb{S}^{d}$.

Kuijlaars and Saff \cite{KuijlaarsSaff:1998Asymptotics} proved that for $d\geq2$ and $s>d$, there exist constants $C_{d,s}^{(1)}, C_{d,s}^{(2)}>0$, such that
\begin{equation}\label{KuijlaarsSaff}
C_{d,s}^{(1)}N^{1+\frac{s}{d}}\leq \mathcal{E}_{s}^{(d)}(N)\leq C_{d,s}^{(2)}N^{1+\frac{s}{d}}.
\end{equation}

Also in  \cite{KuijlaarsSaff:1998Asymptotics} it was showed that for $s=d$ the following formula holds
\begin{equation}\label{KuijlaarsSaff_s=d}
\lim\limits_{N\rightarrow\infty}(N^{2}\log N)^{-1} \mathcal{E}_{s}^{(d)}(N)=\frac{1}{2d}\frac{\Gamma(\frac{d+1}{2})}{\Gamma(\frac{d}{2})\Gamma(\frac{1}{2})}.
\end{equation}

 We show that for every well-separated sequence of
 $N$-point spherical $t$-designs on 
   $\mathbb{S}^{d}$, $d\geq2$,   with $N\asymp t^{d}$ the following relations are true:
 \begin{equation*}
E_{s}^{(d)}(X_{N})\ll N^{1+\frac{s}{d}}, \ s>d
\end{equation*}
and
 \begin{equation*}
\lim\limits_{N\rightarrow\infty}\frac{E_{s}^{(d)}(X_{N})}{N^{2}\log N}=\frac{1}{2d\sqrt{\pi}}\frac{\Gamma(\frac{d}{2}+\frac{1}{2})}{\Gamma(\frac{d}{2})}, \ s=d.
\end{equation*}
Here and further we use the Vinogradov notation $a_{n}\ll b_{n}$  to mean that there exists positive constant $C$  independent of $n$ such that 
 $a_{n}\leq  C b_{n}$  for all $n$.

First, we observe, that since $\mathcal{E}_{s}^{(d)}(N)\leq {E}_{s}^{(d)}(X_N)$ for any $N$-point set, the lower bound in (\ref{KuijlaarsSaff}) provides the lower bound for the $s$-energy of any $N$-point set.
So, asymptotically for Riesz $s$-energy, $s\geq d$, well-separated spherical $t$-designs 
 are as good  as point sets which  minimize the $s$-energy.

This paper is organised as follows:  Section~\ref{prelim} provides basic notations and necessary background for Jacobi polynomials, Section~\ref{mainResults} contains formulation of main results and proofs of theorems.

\section{Preliminaries}
\label{prelim}

In this paper we use the Pochhammer symbol $(a)_{n}$, where
$n\in \mathbb{N}_{0}$ and $a\in \mathbb{R}$, defined by
\begin{equation*}
  (a)_{0}:=1, \quad (a)_{n}:=a(a+1)\ldots(a+n-1)\quad \mathrm{for} \quad
  n\in \mathbb{N},
\end{equation*}
which can be  written in the terms of the gamma function $\Gamma(z)$ by means of
\begin{equation}\label{Pochhammer}
 (a)_{\ell}=\frac{\Gamma(\ell+a)}{\Gamma(a)}.
 \end{equation}
 For fixed $a,b$ the following asymptotic equality is true
 \begin{equation}\label{gamma}
 \frac{\Gamma(n+a)}{\Gamma(n+b)}= n^{a-b}\Big(1+\mathcal{O}\Big(\frac{1}{n}\Big) \Big) \ \ \mathrm{as} \ \ n\rightarrow \infty.
\end{equation}

For any integrable function $f: [-1, 1]\rightarrow \mathbb{R}$ (see, e.g.,
\cite{Mueller1966:spherical_harmonics}) we have
\begin{equation}\label{a1}
 \int\limits_{\mathbb{S}^{d}}f(\langle\mathbf{x},\mathbf{y}\rangle)d\sigma_{d}(\mathbf{x})=\frac{\Gamma(\frac{d+1}{2})}{\sqrt{\pi}\Gamma(\frac{d}{2})}\int\limits_{-1}^{1}f(t)(1-t^{2})^{\frac{d}{2}-1}dt \quad \forall \mathbf{y}\in \mathbb{S}^{d}.
\end{equation}

The Jacobi polynomials $P_{\ell}^{(\alpha,\beta)}(x)$ are the polynomials
orthogonal over the interval $[-1,1]$ with the weight function
$w_{\alpha,\beta}(x)=(1-x)^{\alpha}(1+x)^{\beta}$ and normalised by the
relation
\begin{equation}\label{JacobiMax}
  P_{\ell}^{(\alpha,\beta)}(1)=\binom {\ell+\alpha}\ell=
  \frac{(1+\alpha)_{\ell}}{\ell!}=\frac{1}{\Gamma(1+\alpha)}\ell^{\alpha}\left(1+\mathcal{O}\Big(\frac{1}{\ell}\Big)\right),
  \quad \alpha,\beta>-1.
 \end{equation}
 (see, e.g., \cite[(5.2.1)]{Magnus-Oberhettinger-Soni1966:formulas_theorems}).
 
We will also use formula
\begin{equation}\label{minusArgum}
  P_{\ell}^{(\alpha,\beta)}(-x)=(-1)^{l}
  P_{\ell}^{(\alpha,\beta)}(x)
 \end{equation}
 and the connection coefficient formula (see, e.g.,  Theorem 7.1.4 from \cite{SpecialFunctions})
\begin{equation}\label{JacobiConectForm}
P_{m}^{(\gamma,\gamma)}(x)=\frac{(\gamma+1)_{m}}{(2\gamma+1)_{m}}
\sum\limits_{k=0}^{[\frac{m}{2}]}\frac{(2\alpha+1)_{m-2k}}{(\alpha+1)_{m-2k}}
  \frac{(\gamma+\frac{1}{2})_{m-k}(\alpha+\frac{3}{2})_{m-2k}(\gamma-\alpha)_{k}}{(\alpha+\frac{3}{2})_{m-k}(\alpha+\frac{1}{2})_{m-2k}k!}P_{m-2k}^{\alpha,\alpha}(x).
 \end{equation}

 For fixed ${\alpha, \beta>-1}$ and ${0< \theta<\pi}$, the following relation
 gives an asymptotic approximation for $\ell\rightarrow\infty$ (see,
 e.g., \cite[Theorem 8.21.13]{Szegoe1975:orthogonal_polynomials})
\begin{multline*}
  P_{\ell}^{(\alpha,\beta)}(\cos \theta)=\frac{1}{\sqrt{\pi}}\ell^{-1/2}
  \Big(\sin\frac{\theta}{2}\Big)^{-\alpha-1/2}
  \Big(\cos\frac{\theta}{2}\Big)^{-\beta-1/2}\\
  \times\Big\{\cos \Big(\Big(\ell+\frac{\alpha+\beta+1}{2}\Big)\theta-
  \frac{2\alpha+1}{4}\pi\Big)+\mathcal{O}(\ell\sin\theta)^{-1}\Big\}.
\end{multline*}
Thus, for
$c_{\alpha,\beta}\ell^{-1}\leq\theta\leq\pi-c_{\alpha,\beta}\ell^{-1}$ the last
asymptotic equality yields
 \begin{equation}\label{JacobiIneq}
   |P_{\ell}^{(\alpha,\beta)}(\cos \theta)|\leq \tilde{c}_{\alpha,\beta}
   \ell^{-1/2}(\sin\theta)^{-\alpha-1/2}+
 \tilde{c}_{\alpha,\beta}\ell^{-3/2}(\sin\theta)^{-\alpha-3/2}, \quad\alpha\geq\beta.
 \end{equation}

The following differentiation formula holds
\begin{equation}\label{JacobiDifferen}
   \frac{d}{dx}P_{n}^{(\alpha,\beta)}(x)=\frac{\alpha+\beta+n+1}{2}P_{n-1}^{(\alpha+1,\beta+1)}(x).
 \end{equation}

If $\lambda>s-1$, $s\geq d$, then taking into account formula  \cite[(5.3.4)]{Magnus-Oberhettinger-Soni1966:formulas_theorems}) and the fact  that 
the Gegenbauer polynomials  are a special case of the Jacobi polynomials $P_{n}^{(\alpha,\beta)}(x)$ (see, e.g.,
 \cite[(5.3.1)]{Magnus-Oberhettinger-Soni1966:formulas_theorems}),
we have that  for ${-1<x<1}$ the following expansion holds 
\begin{multline}\label{expansionGegenbauer1}
(1-x)^{-\frac{s}{2}}=2^{2\lambda-\frac{s}{2}}\pi^{-\frac{1}{2}}\Gamma(\lambda)\Gamma\Big(\lambda-\frac{s}{2}+\frac{1}{2}\Big) \\
\times\sum\limits_{n=0}^{\infty}\frac{(n+\lambda)(\frac{s}{2})_{n}}{\Gamma(n+2\lambda-\frac{s}{2}+1)}\frac{(2\lambda)_{n}}{(\lambda+\frac{1}{2})_{n}}P_{n}^{(\lambda-\frac{1}{2},\lambda-\frac{1}{2})}(x).
\end{multline}

\section{Main results}
\label{mainResults}

By a spherical cap $S(\mathbf{x}; \varphi)$ of centre $\mathbf{x}$ and angular radius
 $\varphi$ we mean
 \begin{equation*}
 S(\mathbf{x}; \varphi):=\big\{\mathbf{y}\in \mathbb{S}^{d} \big| \langle\mathbf{x},\mathbf{y}\rangle\geq \cos\varphi \big\}.
 \end{equation*}
The normalised surface area of a spherical cap is given by 
\begin{equation}\label{capArea}
|S(\mathbf{x}; \varphi)|=\frac{\Gamma((d+1)/2)}{\sqrt{\pi}\Gamma(d/2)}
\int\limits_{\cos\varphi}^{1}(1-t^{2})^{\frac{d}{2}-1}dt
\asymp(1-\cos\varphi)^{\frac{d}{2}} \quad\text{as } \varphi\rightarrow 0.
\end{equation}

If for sequence $(X_{N})_{N}$ condition (\ref{wellSeparat}) holds, then any  spherical cap $S(\mathbf{x}; \alpha_{N})$, $\mathbf{x}\in \mathbb{S}^{d}$, where
\begin{equation}\label{alphaN}
\alpha_{N}:=\arccos\Big(1-\frac{c^{2}_{1}}{8N^{\frac{2}{d}}}\Big),
\end{equation}
 contains at
 most one point of the set $(X_{N})_{N}$. 

From the elementary estimates
\begin{equation}\label{sinIneq}
\sin\theta\leq \theta\leq \frac{\pi}{2}\sin\theta, \quad 0\leq\theta\leq \frac{\pi}{2},
\end{equation}
we obtain
\begin{equation}\label{alphaEstim}
\Big(1-\frac{c^{2}_{1}}{16N^{\frac{2}{d}}}\Big)^{\frac{1}{2}}\frac{c_{1}}{2N^{\frac{1}{d}}}\leq\alpha_{N}\leq 
\frac{\pi}{4}\Big(1-\frac{c^{2}_{1}}{16N^{\frac{2}{d}}}\Big)^{\frac{1}{2}}\frac{c_{1}}{N^{\frac{1}{d}}}.
 \end{equation} 

The following two theorems are the main result of this paper. 
 
\begin{theorem}\label{theoremLogar} Let $d\geq2$ be fixed, 
  $(X_{N(t)})_t$ be a sequence of well-separated spherical $t$-designs on
   $\mathbb{S}^{d}$ and $N(t)\asymp t^{d}$.
Then for the logarithmic energy $E_{log}^{(d)}(X_{N})$ the following estimate holds
 \begin{equation}\label{theoremLog}
 E_{log}^{(d)}(X_{N})=N^{2}\int\limits_{\mathbb{S}^{d}}\int\limits_{\mathbb{S}^{d}}\log\frac{1}{|\mathbf{x}-\mathbf{y}|}d\sigma_{d}(\mathbf{x})d\sigma_{d}(\mathbf{y})  -\frac{1}{d}N\log N +\mathcal{O}(N).
\end{equation}

\end{theorem} 
 
\begin{theorem}
\label{theoremRisz} Let $d\geq2$ be fixed, and
  $(X_{N(t)})_t$ be a sequence of well-separated spherical $t$-designs on
   $\mathbb{S}^{d}$ and $N(t)\asymp t^{d}$.
Then for $s>d$ the $s$-energy $E_{s}^{(d)}(X_{N})$ satisfies the estimate 
 \begin{equation}\label{theorems>d}
 E_{s}^{(d)}(X_{N})\ll N^{1+\frac{s}{d}}, 
\end{equation}
and for $s=d$, the $s$-energy $E_{s}^{(d)}(X_{N})$ satisfies following estimates
 \begin{equation}\label{theorems=d}
 E_{s}^{(d)}(X_{N})=\frac{1}{2\sqrt{\pi}}\frac{\Gamma(\frac{d}{2}+\frac{1}{2})}{\Gamma(\frac{d}{2})}\sum\limits_{n=1}^{[\frac{t}{2}]}n^{-1}N^{2}+\mathcal{O}(N^{2})
\end{equation}
and 
\begin{equation}\label{theorems=dLim}
 \lim\limits_{N\rightarrow\infty}\frac{E_{s}^{(d)}(X_{N})}{N^{2}\log N}=\frac{1}{2d\sqrt{\pi}}\frac{\Gamma(\frac{d}{2}+\frac{1}{2})}{\Gamma(\frac{d}{2})}.
\end{equation}
\end{theorem} 
 

 \subsection{Proof of Theorem  \ref{theoremLogar}} 

For each $i\in\{1,\ldots,N\}$ we divide the sphere $\mathbb{S}^{d}$ into an
  upper hemisphere $H_{i}^{+}$ with 'north pole' $\mathbf{x}_{i}$ and a lower
  hemisphere $H_{i}^{-}$:
\begin{equation*}
H_{i}^{+}:=\Big\{\mathbf{x}\in\mathbb{S}^{d}\Big|\langle\mathbf{x}_{i},\mathbf{x}\rangle\geq0 \Big\},
\end{equation*}
\begin{equation*}
H_{i}^{-}:=\mathbb{S}^{d}\setminus H_{i}^{+}.
\end{equation*}

Noting that
\begin{equation}\label{distance}
|\mathbf{x}_{i}-\mathbf{x}_{j}|^{-1}=\frac{1}{\sqrt{2}}(1-\langle\textbf{x}_{i},\mathbf{x}_{j}\rangle)^{-\frac{1}{2}},
\end{equation}
the logarithmic energy can be written in the form
\begin{equation}\label{formConnect}
 E_{log}^{(d)}(X_{N})=
  {\mathop{\sum}\limits_{i,j=1,\atop i\neq j}^{N}}
  \log \frac{1}{|\mathbf{x}_{i}-\mathbf{x}_{j}|}=
  \frac{1}{2}{\mathop{\sum}\limits_{i,j=1,\atop i\neq j}^{N}}\left(\log \frac{1}{1-\langle\mathbf{x}_{i},\mathbf{x}_{j}\rangle}-\log 2 \right).
\end{equation}

Let $\lambda>d+1$. The, putting $s=2$ in (\ref{expansionGegenbauer1}), we get
\begin{align}\label{for3}
(1-x)^{-1}=2^{2\lambda-1}\pi^{-\frac{1}{2}}\Gamma(\lambda)\Gamma\Big(\lambda-\frac{1}{2}\Big)\sum\limits_{n=0}^{\infty}\frac{(n+\lambda)\Gamma(n+1)}{\Gamma(n+2\lambda)}
\frac{(2\lambda)_{n}}{(\lambda+\frac{1}{2})_{n}}P_{n}^{(\lambda-\frac{1}{2}, \lambda-\frac{1}{2})}(x).
\end{align}

Formula (\ref{JacobiDifferen}) implies, that
\begin{equation}\label{integration_formula}
\int P_{n}^{(\lambda-\frac{1}{2}, \lambda-\frac{1}{2})}(x)dx=\frac{2}{n+2\lambda-1}
P_{n+1}^{(\lambda-\frac{3}{2}, \lambda-\frac{3}{2})}(x).
\end{equation}
Integrating from $0$ to $x$, we have
\begin{align}
\log \frac{1}{1-x}&=2^{2\lambda}\pi^{-\frac{1}{2}}\Gamma(\lambda)\Gamma\Big(\lambda-\frac{1}{2}\Big) \notag \\
& \times\sum\limits_{n=0}^{\infty}\frac{(n+\lambda)\Gamma(n+1)}{(n+2\lambda-1)\Gamma(n+2\lambda)}
\frac{(2\lambda)_{n}}{(\lambda+\frac{1}{2})_{n}}
\big(P_{n}^{(\lambda-\frac{1}{2}, \lambda-\frac{1}{2})}(x)-P_{n}^{(\lambda-\frac{1}{2}, \lambda-\frac{1}{2})}(0)\big). \label{for4}
\end{align}

 We split the $\log$-energy into two parts
 \begin{equation}\label{1splitLog}
E_{log}^{(d)}(X_{N})=\sum\limits_{j=1}^{N}
  {\mathop{\sum}\limits_{i=1,\atop \mathbf{x}_{i}\in H^{\pm}_{i}\setminus S(\pm \mathbf{x}_{j};\alpha_{N})}^{N}} \log \frac{1}{|\mathbf{x}_{i}-\mathbf{x}_{j}|}
  +
  \sum\limits_{j=1}^{N}
  {\mathop{\sum}\limits_{i=1,\atop \mathbf{x}_{i}\in  S(- \mathbf{x}_{j};\alpha_{N})}^{N}}
   \log \frac{1}{|\mathbf{x}_{i}-\mathbf{x}_{j}|}.
\end{equation}

From (\ref{wellSeparat}) and the fact the spherical cap $S(- \mathbf{x}_{j};\alpha_{N})$ contains at most one point of $X_{N}$, the second term in (\ref{1splitLog}), where the scalar product is close to $-1$, can be bounded from above by
\begin{equation}\label{estim1split}
 \sum\limits_{j=1}^{N}
  {\mathop{\sum}\limits_{i=1,\atop \mathbf{x}_{i}\in  S(- \mathbf{x}_{j};\alpha_{N})}^{N}}\log \frac{1}{|\mathbf{x}_{i}-\mathbf{x}_{j}|}=\mathcal{O}(N).
\end{equation}

Taking into account (\ref{formConnect}), (\ref{for4})--(\ref{estim1split}), we deduce
\begin{equation}\label{formSubst}
 E_{log}^{(d)}(X_{N})=
  \frac{1}{2}E_{H_{\log,t}}(X_{N})+\frac{1}{2}E_{R_{\log,t}}(X_{N})-
  \frac{1}{2}N^{2}\log 2 +\mathcal{O}(N),
\end{equation}
where
\begin{align}
&H_{\log,t}(x)=H_{\log,t}(d,\lambda,x) \notag\\
&:=-
2^{2\lambda}\pi^{-\frac{1}{2}}\Gamma(\lambda)\Gamma\Big(\lambda-\frac{1}{2}\Big)\sum\limits_{n=0}^{\infty}\frac{(n+\lambda)\Gamma(n+1)}{(n+2\lambda-1)\Gamma(n+2\lambda)}
\frac{(2\lambda)_{n}}{(\lambda+\frac{1}{2})_{n}}
P_{n+1}^{(\lambda-\frac{3}{2},\lambda-\frac{3}{2})}(0) \notag\\
&+
2^{2\lambda}\pi^{-\frac{1}{2}}\Gamma(\lambda)\Gamma\Big(\lambda-\frac{1}{2}\Big)\sum\limits_{n=0}^{t-1}\frac{(n+\lambda)\Gamma(n+1)}{(n+2\lambda-1)\Gamma(n+2\lambda)}
\frac{(2\lambda)_{n}}{(\lambda+\frac{1}{2})_{n}}
P_{n+1}^{(\lambda-\frac{3}{2},\lambda-\frac{3}{2})}(x),\label{hLog_tDefinition}
\end{align}
\begin{align}
&R_{\log,t}(x)=R_{\log,t}(d,\lambda,x)\notag\\
&
:= 2^{2\lambda}\pi^{-\frac{1}{2}}\Gamma(\lambda)\Gamma\Big(\lambda-\frac{1}{2}\Big)\sum\limits_{n=t+1}^{\infty}\frac{(n+\lambda)\Gamma(n+1)}{(n+2\lambda-1)\Gamma(n+2\lambda)}
\frac{(2\lambda)_{n}}{(\lambda+\frac{1}{2})_{n}}
P_{n+1}^{(\lambda-\frac{3}{2},\lambda-\frac{3}{2})}(x)\label{rLog_tDefinition}
\end{align}
and 
\begin{equation}\label{energyParts}
E_{U}(X):=\sum\limits_{j=1}^{N}
  {\mathop{\sum}\limits_{i=1,\atop \textbf{x}_{i}\in H^{\pm}_{i}\setminus S(\pm x_{j};\alpha_{N})}^{N}}U(\langle\textbf{x}_{i},\mathbf{x}_{j}\rangle).
\end{equation}

Let us show that
\begin{equation}\label{estimR}
E_{R_{\log,t}}(X_{N})=\mathcal{O}(N).
\end{equation}

Applying (\ref{Pochhammer}), (\ref{gamma}) and (\ref{JacobiIneq}) to (\ref{rLog_tDefinition}),
we have
\begin{align}
|R_{\log,t}(\cos\theta)| \ll&
\sum\limits_{n=t}^{\infty}n^{-\lambda+\frac{1}{2}}
|P_{n+1}^{\lambda-\frac{3}{2},\lambda-\frac{3}{2}}(\cos\theta)| \notag \\
\ll&\sum\limits_{n=t}^{\infty}n^{-\lambda+\frac{1}{2}}\left( n^{-\frac{1}{2}}(\sin\theta)^{-\lambda+1} +n^{-\frac{3}{2}}(\sin\theta)^{-\lambda}\right) \notag \\
\ll &
(\sin\theta)^{-\lambda+1}t^{-\lambda+1}+(\sin\theta)^{-\lambda}t^{-\lambda}.
\label{estim1}
\end{align}
 From \cite[(3.30) and (3.33)]{Brauchart-Hesse2007:numerical_integration}, it
 follows that
 \begin{multline} \label{BrauchartHesse}
  \sum\limits_{j=1}^{N}{\mathop{\sum}\limits^{N}_{
 i=1,\atop \mathbf{x}_{i}\in H_{j}^{\pm}\setminus
S(\pm\mathbf{x}_{j}; \frac{c}{n})}}
  (\sin\theta_{ij}^{\pm})^{-\frac{d}{2}+\frac{1}{2}-k-L}
\\
 \ll N^{2}(1+n^{L+k-(d+1)/2}), \quad k=0,1,\ldots \quad \text{for }L>\frac{d+1}{2}.
 \end{multline}

Estimates (\ref{alphaEstim}) and (\ref{BrauchartHesse}) imply
\begin{align}
& E_{R_{\log,t}}(X_{N})\ll t^{-\lambda+1}\sum\limits_{j=1}^{N}
  {\mathop{\sum}\limits_{i=1,\atop \mathbf{x}_{i}\in H^{\pm}_{i}\setminus S(\pm \mathbf{x}_{j};\alpha_{N})}^{N}}   (\sin\theta^{\pm})^{-\lambda+1}+t^{-\lambda}\sum\limits_{j=1}^{N}
  {\mathop{\sum}\limits_{i=1,\atop \mathbf{x}_{i}\in H^{\pm}_{i}\setminus S(\pm \mathbf{x}_{j};\alpha_{N})}^{N}}   (\sin\theta^{\pm})^{-\lambda}
\notag\\
&\ll N^{2}t^{-d}\ll N, \ \ \lambda>d+1.
\label{estim2}
\end{align}
This proves (\ref{estimR}).

Now let us find the estimate for $E_{H_{\log,t}}(X_{N})$.
The polynomial $H_{log,t}$ is a spherical polynomial of degree $t$ and $X_{N}$ is a spherical $t$-design. That is why,  an equal weight integration rule with nodes from $X_{N}$ integrates $H_{log,t}$ exactly, and 
\begin{align}
& E_{H_{\log,t}}(X_{N})=\sum\limits_{j=1}^{N}
  {\mathop{\sum}\limits_{i=1,\atop \mathbf{x}_{i}\in H^{\pm}_{i}\setminus S(\pm \mathbf{x}_{j};\alpha_{N})}^{N}} H_{\log,t}(\langle\mathbf{x}_{i},\mathbf{x}_{j}\rangle)
\notag\\
&=
N^{2}\int\limits_{\mathbb{S}^{d}}H_{\log,t}(\langle\mathbf{x},\mathbf{y}\rangle)d\sigma_{d}(\mathbf{x})
-N H_{\log,t}(1)-  \sum\limits_{j=1}^{N}
  {\mathop{\sum}\limits_{i=1,\atop \mathbf{x}_{i}\in S(- \mathbf{x}_{j};\alpha_{N})}^{N}} H_{\log,t}(\langle\mathbf{x}_{i},\mathbf{x}_{j}\rangle), \ \mathbf{y}\in \mathbb{S}^{d} .
\label{estim3}
\end{align}

Let $b_{0}\in\mathbb{R}_{+}$ is such, that for $\beta_{N}:=\arccos(1-b_{0}N^{-\frac{2}{d}})$ the following relation holds
\begin{equation}\label{betaN}
\int\limits_{S(\mathbf{y};\beta_{N})} d\sigma_{d}(\mathbf{x})=\frac{\Gamma(\frac{d+1}{2})}{\sqrt{\pi}\Gamma(\frac{d}{2})}\int\limits_{1-b_{0}N^{-\frac{2}{d}}}^{1}(1-x^{2})^{\frac{d}{2}-1}dx=\frac{1}{N}, \ \ \ \mathbf{y}\in \mathbb{S}^{d}.
\end{equation}

It is clear, that
\begin{equation}\label{betaNEstim}
\beta_{N}\asymp N^{-\frac{1}{d}}.
\end{equation}

Then
\begin{equation}\label{equal1}
 E_{H_{\log,t}}(X_{N})=
 N^{2}\int\limits_{\mathbb{S}^{d}}\log\frac{1}{1-\langle\mathbf{x},\mathbf{y}\rangle}d\sigma_{d}(\mathbf{x})+Q_{t}(X_{N}),
\end{equation}
where
\begin{align}
& 
Q_{t}(X_{N})=Q_{t}(d,X_{N}):= -N^{2}\int\limits_{S(\pm\mathbf{y};\beta_{N})}\log\frac{1}{1-\langle\mathbf{x},\mathbf{y}\rangle}d\sigma_{d}(\mathbf{x})
\notag \\
&-N^{2}\int\limits_{\mathbb{S}^{d}\setminus S(\pm\mathbf{y};\beta_{N})}R_{\log,t}(\langle\mathbf{x},\mathbf{y}\rangle) d\sigma_{d}(\mathbf{x})+ 
N^{2}\int\limits_{ S(\pm\mathbf{y};\beta_{N})}H_{\log,t}(\langle\mathbf{x},\mathbf{y}\rangle) d\sigma_{d}(\mathbf{x})
\notag\\
&
-N H_{\log,t}(1)-  \sum\limits_{j=1}^{N}
  {\mathop{\sum}\limits_{i=1,\atop \mathbf{x}_{i}\in S(- \mathbf{x}_{j};\alpha_{N})}^{N}} H_{\log,t}(\langle\mathbf{x}_{i},\mathbf{x}_{j}\rangle), \ \ \mathbf{y}\in \mathbb{S}^{d}.
\label{estim4}
\end{align}
Now we shall prove that 
\begin{equation}\label{estimQ}
Q_{t}(X_{N})=-N^{2}\int\limits_{S(\mathbf{y};\beta_{N})}\log\frac{1}{1-\langle\mathbf{x},\mathbf{y}\rangle}d\sigma_{d}(\mathbf{x})+\mathcal{O}(N), \ \ \mathbf{y}\in \mathbb{S}^{d}.
\end{equation}

Using (\ref{a1}), (\ref{estim1}) and (\ref{betaNEstim}), we get
\begin{align}
& 
N^{2}\left| \int\limits_{\mathbb{S}^{d}\setminus S(\pm\mathbf{y};\beta_{N})}R_{\log,t}(\langle\mathbf{x},\mathbf{y}\rangle) d\sigma_{d}(\mathbf{x}) \right| \ll
N^{2}\int\limits_{-1+b_{0}N^{-\frac{2}{d}}}^{1-b_{0}N^{-\frac{2}{d}}} |R_{\log,t}(x)|(1-x^{2})^{\frac{d}{2}-1}dx \notag \\
& \ll N^{2}\int\limits_{-1+b_{0}N^{-\frac{2}{d}}}^{1-b_{0}N^{-\frac{2}{d}}} 
\left(t^{-\lambda+1}(\sqrt{1-x^{2}})^{-\lambda+1}+t^{-\lambda}(\sqrt{1-x^{2}})^{-\lambda} \right)
(1-x^{2})^{\frac{d}{2}-1}dx \notag \\
& = 2N^{2}\int\limits_{\beta_{N}}^{\frac{\pi}{2}}\left(t^{-\lambda+1}(\sin y)^{-\lambda+1}+
t^{-\lambda}(\sin y)^{-\lambda}\right)(\sin y)^{d-1} dy \notag \\
& 
\ll N^{2}\int\limits_{\beta_{N}}^{\frac{\pi}{2}}\left(t^{-\lambda+1}y^{-\lambda+d}+
t^{-\lambda}y^{-\lambda+d-1}\right) dy \ll N.
\label{estim5}
\end{align}

From the definition of $\beta_{n}$ it is easy to see, that
\begin{equation}\label{estim6}
\left| N^{2}\int\limits_{S(-\mathbf{y};\beta_{N})}\log\frac{1}{1-\langle\mathbf{x},\mathbf{y}\rangle}d\sigma_{d}(\mathbf{x})\right|\ll N^{2}|S(-\mathbf{y};\beta_{N})|\ll N, \ \mathbf{y}\in \mathbb{S}^{d}.
\end{equation}

According to the definition of $\beta_N$ (\ref{betaN}) we deduce
\begin{align}
& 
\left| N^{2}\int\limits_{ S(\mathbf{y};\beta_{N})}H_{\log,t}(\langle\mathbf{x},\mathbf{y}\rangle) d\sigma_{d}(\mathbf{x}) - N H_{\log,t}(1)\right| \notag \\
&=\left| N^{2}\frac{\Gamma(\frac{d+1}{2})}{\sqrt{\pi}\Gamma(\frac{d}{2})}\int\limits_{1-b_{0}N^{-\frac{2}{d}}}^{1} (H_{\log,t}(x)-H_{\log,t}(1)) (1-x^{2})^{\frac{d}{2}-1}dx \right| \notag \\
&\ll N \max\limits_{x\in[1-b_{0}N^{-\frac{2}{d}}, 1]}
\left( H_{\log,t}(1)-H_{\log,t}(x) \right)\ll N^{1-\frac{2}{d}}|H^{'}_{\log,t}(1)|.
\label{estim7}
\end{align}

Formulas (\ref{JacobiMax}), (\ref{JacobiDifferen}) and (\ref{hLog_tDefinition}) imply
\begin{align}
& H^{'}_{\log,t}(1)=
2^{2\lambda-1}\pi^{-\frac{1}{2}}\Gamma(\lambda)\Gamma\Big(\lambda-\frac{1}{2}\Big)
 \notag \\
&\times
\sum\limits_{n=0}^{t-1}\frac{(n+\lambda)\Gamma(n+1)}{\Gamma(n+2\lambda)}\frac{(2\lambda)_{n}}{(\lambda+\frac{1}{2})_{n}}
P_{n}^{(\lambda-\frac{1}{2},\lambda-\frac{1}{2})}(1) \ll t^{2}\ll N^{\frac{2}{d}}.
\label{hLog_tDerivative}
\end{align}

From (\ref{minusArgum}), (\ref{JacobiIneq}) and (\ref{hLog_tDefinition}) it follows that
\begin{align}
&|H_{\log,t}(-1)|\ll
\sum\limits_{n=0}^{\infty}\frac{(n+\lambda)\Gamma(n+1)}{(n+2\lambda-1)\Gamma(n+2\lambda)}
\frac{(2\lambda)_{n}}{(\lambda+\frac{1}{2})_{n}}
\frac{1}{\sqrt{n+1}} \notag\\
&+
\left|\sum\limits_{n=0}^{t-1}(-1)^{n+1}\frac{(n+\lambda)\Gamma(n+1)}{(n+2\lambda-1)\Gamma(n+2\lambda)}
\frac{(2\lambda)_{n}}{(\lambda+\frac{1}{2})_{n}}
P_{n+1}^{(\lambda-\frac{3}{2},\lambda-\frac{3}{2})}(1)\right|.\label{hLog_minus1}
\end{align}

Relations (\ref{gamma}) and (\ref{JacobiMax}) allow us to write 
\begin{align}
&\frac{(n+\lambda)\Gamma(n+1)}{(n+2\lambda-1)\Gamma(n+2\lambda)}
\frac{(2\lambda)_{n}}{(\lambda+\frac{1}{2})_{n}}
P_{n+1}^{(\lambda-\frac{3}{2},\lambda-\frac{3}{2})}(1) \notag\\
&=
\frac{\Gamma(\lambda+\frac{1}{2})}{\Gamma(2\lambda)}\frac{1}{n+2\lambda-1}\Big(1+\frac{\lambda-1}{n+1}\Big)\Big(1+\mathcal{O}\Big(\frac{1}{n+1}\Big)\Big)^{2}. \label{internEst}
\end{align}

Hence, (\ref{hLog_minus1}) and (\ref{internEst}) enable us to obtain
\begin{equation}\label{hLog_minus2}
|H_{\log,t}(-1)|=\mathcal{O}(N).
\end{equation}

Using (\ref{a1}), (\ref{hLog_tDerivative}) and (\ref{hLog_minus2}), we deduce
\begin{align}
& 
\left| N^{2}\int\limits_{ S(-\mathbf{y};\beta_{N})}H_{\log,t}(\langle\mathbf{x},\mathbf{y}\rangle) d\sigma_{d}(\mathbf{x}) \right| \notag \\
&=\left| N^{2}\frac{\Gamma(\frac{d+1}{2})}{\sqrt{\pi}\Gamma(\frac{d}{2})}\int\limits_{1-b_{0}N^{-\frac{2}{d}}}^{1} (H_{\log,t}(-x)-H_{\log,t}(-1)) (1-x^{2})^{\frac{d}{2}-1}dx+N H_{\log,t}(-1)\right| \notag \\
& \ll N^{1-\frac{2}{d}}|H^{'}_{\log,t}(1)|+N\ll N.
\label{estim8}
\end{align}

Applying (\ref{hLog_minus2}), we have
\begin{align}
&\left| \sum\limits_{j=1}^{N}
  {\mathop{\sum}\limits_{i=1,\atop \mathbf{x}_{i}\in S(- \mathbf{x}_{j};\alpha_{N})}^{N}} H_{\log,t}(\langle\mathbf{x}_{i},\mathbf{x}_{j}\rangle)   \right|\ll N|H_{\log,t}(\xi)| \notag \\
  &=
  N|H_{\log,t}(\xi)-H_{\log,t}(-1)+H_{\log,t}(-1)| =\mathcal{O}(N), \label{estim9}
\end{align}
where $\xi\in[-1, -1+b_{0}N^{-\frac{2}{d}}]$.

Relations (\ref{estim5})-(\ref{hLog_tDerivative}), (\ref{estim8}) and (\ref{estim9}) prove (\ref{estimQ}).

Integrating by parts, we obtain
\begin{align}
& N^{2}\int\limits_{S(\mathbf{y};\beta_{N})}\log\frac{1}{1-\langle\mathbf{x},\mathbf{y}\rangle}d\sigma_{d}(\mathbf{x})=
\frac{\Gamma(\frac{d+1}{2})}{\sqrt{\pi}\Gamma(\frac{d}{2})}\int\limits_{1-b_{0}N^{-\frac{2}{d}}}^{1} \log\frac{1}{1-x} (1-x^{2})^{\frac{d}{2}-1}dx
 \notag \\
  &
  = \frac{\Gamma(\frac{d+1}{2})}{\sqrt{\pi}\Gamma(\frac{d}{2})}\int\limits_{1-b_{0}N^{-\frac{2}{d}}}^{1} \log\frac{1}{1-x} d\left(-\int\limits_{x}^{1}(1-t^{2})^{\frac{d}{2}-1}dt\right)dx \notag \\
  & = N\log\Big(\frac{N^{\frac{2}{d}}}{b_{0}} \Big)+
  N^{2}\frac{\Gamma(\frac{d+1}{2})}{\sqrt{\pi}\Gamma(\frac{d}{2})}\int\limits_{1-b_{0}N^{-\frac{2}{d}}}^{1} \frac{1}{1-x} \int\limits_{x}^{1}(1-t^{2})^{\frac{d}{2}-1}dt dx \notag \\
  &=
  \frac{2}{d}N\log N+\mathcal{O}(N).
  \label{estim10}
\end{align}

So, combining (\ref{formSubst}), (\ref{estimR}),  (\ref{equal1}), (\ref{estimQ})  and (\ref{estim10}), we get
\begin{align}
 &E_{log}^{(d)}(X_{N})=
  \frac{1}{2}N^{2}\int\limits_{\mathbb{S}^{d}}\int\limits_{\mathbb{S}^{d}}\log\frac{1}{1-\langle\mathbf{x},\mathbf{y}\rangle}d\sigma_{d}(\mathbf{x})d\sigma_{d}(\mathbf{y})  -\frac{1}{d}N\log N-
  \frac{1}{2}N^{2}\log 2 +\mathcal{O}(N) \notag \\
  &= N^{2}\int\limits_{\mathbb{S}^{d}}\int\limits_{\mathbb{S}^{d}}
  \log\frac{1}{|\mathbf{x}-\mathbf{y}|}d\sigma_{d}(\mathbf{x})d\sigma_{d}(\mathbf{y})  -\frac{1}{d}N\log N +\mathcal{O}(N)
  .
\end{align}\label{summar}

This implies (\ref{theoremLog}).
Theorem~\ref{theoremLogar} is proved.

 \qed

\subsection{Proof of Theorem  \ref{theoremRisz}}

In the same way as in the case for logarithmic energy, we split the $s$-energy into two parts
\begin{align}\label{1split}
E_{s}^{(d)}(X_{N})&=\frac{1}{2}\sum\limits_{j=1}^{N}
  {\mathop{\sum}\limits_{i=1,\atop \textbf{x}_{i}\in H^{\pm}_{i}\setminus S(\pm \mathbf{x}_{j};\alpha_{N})}^{N}}|\mathbf{x}_{i}-\mathbf{x}_{j}|^{-s}+
  \frac{1}{2}\sum\limits_{j=1}^{N}
  {\mathop{\sum}\limits_{i=1,\atop \mathbf{x}_{i}\in  S(- \mathbf{x}_{j};\alpha_{N})}^{N}}|\mathbf{x}_{i}-\mathbf{x}_{j}|^{-s} \notag \\
  &=\frac{1}{2}\sum\limits_{j=1}^{N}
  {\mathop{\sum}\limits_{i=1,\atop \textbf{x}_{i}\in H^{\pm}_{i}\setminus S(\pm \mathbf{x}_{j};\alpha_{N})}^{N}}|\mathbf{x}_{i}-\mathbf{x}_{j}|^{-s}+\mathcal{O}(N).
\end{align}

Taking into account that the Jacobi  series (\ref{expansionGegenbauer1}) converges uniformly in \linebreak ${\Big[-1+\frac{c^{2}_{1}}{8N^{\frac{2}{d}}},1-\frac{c^{2}_{1}}{8N^{\frac{2}{d}}} \Big]}$,
 for $\lambda>s-1$  we get that 
\begin{align}\label{expansionSubst}
\frac{1}{2}\sum\limits_{j=1}^{N}
  {\mathop{\sum}\limits_{i=1,\atop \mathbf{x}_{i}\in H^{\pm}_{i}\setminus S(\pm \mathbf{x}_{j};\alpha_{N})}^{N}}|\mathbf{x}_{i}-\mathbf{x}_{j}|^{-s}&=
  \frac{1}{2^{1+\frac{s}{2}}}\sum\limits_{j=1}^{N}
  {\mathop{\sum}\limits_{i=1,\atop \mathbf{x}_{i}\in H^{\pm}_{i}\setminus S(\pm \mathbf{x}_{j};\alpha_{N})}^{N}}(1-\langle\textbf{x}_{i},\mathbf{x}_{j}\rangle)^{-\frac{s}{2}} \notag \\
  &=\frac{1}{2}E_{H_{s,t}}(X_{N})+\frac{1}{2}E_{R_{s,t}}(X_{N}),
\end{align}
where
\begin{align}\label{s_tDefinition}
&H_{s,t}(x)=H_{s,t}(d,\lambda,x) \notag \\
&:=2^{2\lambda-s}\pi^{-\frac{1}{2}}\Gamma(\lambda)\Gamma\Big(\!\lambda\!-\!\frac{s}{2}\!+\!\frac{1}{2}\!\Big)\sum\limits_{n=0}^{t}\frac{(n+\lambda)(\frac{s}{2})_{n}}{\Gamma(n+2\lambda-\frac{s}{2}+1 )}\frac{(2\lambda)_{n}}{(\lambda+\frac{1}{2})_{n}}
P_{n}^{(\lambda-\frac{1}{2}, \ \lambda-\frac{1}{2})}(x),
\end{align}

\begin{align}\label{r_tDefinition}
&R_{s,t}(x)=R_{s,t}(d,\lambda,x) \notag \\
&:=2^{2\lambda-s}\pi^{-\frac{1}{2}}\Gamma(\lambda)\Gamma\Big(\!\lambda\!-\!\frac{s}{2}\!+\!\frac{1}{2}\!\Big)\!\sum\limits_{n=t+1}^{\infty}\!\frac{(n+\lambda)(\frac{s}{2})_{n}}{\Gamma(n+2\lambda-\frac{s}{2}+1 )}\frac{(2\lambda)_{n}}{(\lambda+\frac{1}{2})_{n}}
P_{n}^{(\lambda-\frac{1}{2}, \ \lambda-\frac{1}{2})}(x).
\end{align}

Formula (65) from \cite{GrabnerStepanyukJAT} implies
\begin{equation}\label{form1}
E_{R_{s,t}}(X_{N})=\mathcal{O}\Big(N^{1+\frac{s}{d}}\Big).
\end{equation}

Hence, 
\begin{equation}\label{form11}
E_{s}^{(d)}(X_{N})=\frac{1}{2}E_{H_{s,t}}(X_{N})+\mathcal{O}\Big(N^{1+\frac{s}{d}}\Big), \ \ \ \lambda>s-1,
\end{equation}
where we have used formulas (\ref{1split}), (\ref{expansionSubst}) and (\ref{form1}).

The polynom $H_{s,t}$ is a spherical polynomial of degree $t$ and $X_{N}$ is a spherical $t$-design. So,  an equal weight integration rule with nodes from $X_{N}$ integrates $H_{s,t}$ exactly, and 
\begin{align}\label{s_eq1}
\frac{1}{2}E_{H_{s,t}}(X_{N})=& 
\frac{1}{2}\sum\limits_{j=1}^{N}
 \sum\limits_{i=1}^{N}H_{s,t}(\langle\mathbf{x}_{i},\mathbf{x}_{j}\rangle)
-
\frac{1}{2}\sum\limits_{j=1}^{N}
  {\mathop{\sum}\limits_{i=1,\atop \mathbf{x}_{i}\in  S(\pm \mathbf{x}_{j};\alpha_{N})}^{N}}H_{s,t}(\langle\mathbf{x}_{i},\mathbf{x}_{j}\rangle) 
  +
  \mathcal{O}\big(NH_{s,t}(1)\big)
\notag \\
=&\frac{1}{2}N^{2}\int\limits_{\mathbb{S}^{d}}H_{s,t}(\langle\textbf{x},\mathbf{y}\rangle)d\sigma_{d}(\textbf{x})
+\mathcal{O}\big(NH_{s,t}(1)\big), \ \ \mathbf{y}\in\mathbb{S}^{d}
\end{align}

From relations  (\ref{Pochhammer}), (\ref{gamma}), (\ref{JacobiMax}) and (\ref{s_tDefinition}) we obtain
\begin{align}\label{abs_valS_t}
&H_{s,t}(1) \notag
 \\
&= 
2^{2\lambda-s}\pi^{-\frac{1}{2}}\Gamma(\lambda)\Gamma\Big(\!\lambda\!-\!\frac{s}{2}\!+\!\frac{1}{2}\!\Big)
\sum\limits_{n=0}^{t}\frac{(n+\lambda)(\frac{s}{2})_{n}}{\Gamma(n+2\lambda-\frac{s}{2}+1 )}\frac{(2\lambda)_{n}}{(\lambda+\frac{1}{2})_{n}}
P_{n}^{(\lambda-\frac{1}{2}, \ \lambda-\frac{1}{2})}(1) \notag
 \\
&=2^{2\lambda-s}\pi^{-\frac{1}{2}}\Gamma(\lambda)\Gamma\Big(\!\lambda\!-\!\frac{s}{2}\!+\!\frac{1}{2}\!\Big)
\sum\limits_{n=0}^{t}\frac{(n+\lambda)(\frac{s}{2})_{n}}{\Gamma(n+2\lambda-\frac{s}{2}+1 )}\frac{(2\lambda)_{n}}{n!}\ll t^{s}\ll N^{\frac{s}{d}}.
\end{align}

Let now estimate the integral from (\ref{s_eq1}).
Substituting $\gamma=\lambda-\frac{1}{2}$, $\alpha=\frac{d}{2}-1$ in formula (\ref{JacobiConectForm}), we have
\begin{align}\label{JacobiConectForm_dK}
  &P_{n}^{(\lambda-\frac{1}{2},\lambda-\frac{1}{2})}(x) \notag
  \\
&=\frac{(\lambda+\frac{1}{2})_{n}}{(2\lambda)_{n}} \sum\limits_{k=0}^{[\frac{n}{2}]}\frac{(d-1)_{n-2k}}{(\frac{d}{2})_{n-2k}}
  \frac{(\lambda)_{n-k}(\frac{d}{2}+\frac{1}{2})_{n-2k}(\lambda-\frac{d}{2}+\frac{1}{2})_{k}}{(\frac{d}{2}+\frac{1}{2})_{n-k}(\frac{d}{2}-\frac{1}{2})_{n-2k}k!}P_{n-2k}^{\frac{d}{2}-1,\frac{d}{2}-1}(x).
 \end{align}
 
Since
\begin{equation}
 \int\limits_{\mathbb{S}^{d}}P_{n}^{(\frac{d}{2}-1,\frac{d}{2}-1)}(\mathbf{x})d\sigma_{d}(\textbf{x})=0, \ \ \  n\geq 1,
\end{equation} 
 then (\ref{JacobiConectForm_dK}) yields
 \begin{equation}\label{JacobiIntegral00}
  \int\limits_{\mathbb{S}^{d}}P_{n}^{(\lambda-\frac{1}{2}, \lambda-\frac{1}{2})}(\mathbf{x})d\sigma_{d}(\textbf{x})
=\begin{cases}
  0 & \text{if }
  n=2m+1, \\
 \frac{(\lambda+\frac{1}{2})_{n}}{(2\lambda)_{n}}  
\frac{(\lambda)_{\frac{n}{2}}(\lambda-\frac{d}{2}+\frac{1}{2})_{\frac{n}{2}}}{(\frac{d}{2}+\frac{1}{2})_{\frac{n}{2}}(\frac{n}{2})!} 
 & \text{if }  n=2m.
  \end{cases}  
 \end{equation}

 So,
\begin{multline}\label{for2}
\int\limits_{\mathbb{S}^{d}}H_{s,t}(\langle\textbf{x},\mathbf{y}\rangle)d\sigma_{d}(\textbf{x}) \\
=
2^{2\lambda-s}\pi^{-\frac{1}{2}}\Gamma(\lambda)\Gamma \Big(\lambda-\frac{s}{2}+\frac{1}{2}\Big)
\sum\limits_{n=0}^{[\frac{t}{2}]}
\frac{(2n+\lambda)(\frac{s}{2})_{2n}}{\Gamma(2n+2\lambda-\frac{s}{2}+1)}
\frac{(\lambda)_{n}(\lambda-\frac{d}{2}+\frac{1}{2})_{n}}{(\frac{d}{2}+\frac{1}{2})_{n}n!},
\end{multline}
where we have used (\ref{hLog_tDefinition}) and (\ref{JacobiIntegral00}).

Thus, if $s>d$, then
\begin{equation}\label{integrEst_s>d}
\int\limits_{\mathbb{S}^{d}}H_{s,t}(\langle\textbf{x},\mathbf{y}\rangle)d\sigma_{d}(\textbf{x}) \ll t^{s-d} \ll N^{-1+\frac{s}{d}}
\end{equation}
and the
 relations (\ref{s_eq1}), (\ref{abs_valS_t}) and (\ref{integrEst_s>d}) imply
\begin{align}\label{ineq_s_t1}
E_{s}^{(d)}(X_{N})\ll N^{1+s} .
\end{align}
This implies (\ref{theorems>d}).

If $s=d$, then using (\ref{Pochhammer}) and (\ref{gamma}) from (\ref{for2}) we have
\begin{align}\label{integrEst_s=d}
&\int\limits_{\mathbb{S}^{d}}H_{d,t}(\langle\textbf{x},\mathbf{y}\rangle)d\sigma_{d}(\textbf{x}) \notag\\
&=
2^{2\lambda-d}\pi^{-\frac{1}{2}}
\frac{\Gamma (\frac{d}{2}+\frac{1}{2})}{\Gamma (\frac{d}{2})}
\sum\limits_{n=0}^{[\frac{t}{2}]}
\frac{(2n+\lambda)\Gamma(2n+\frac{d}{2})}{\Gamma(2n+2\lambda-\frac{d}{2}+1)}
\frac{\Gamma(n+\lambda)\Gamma(n+\lambda-\frac{d}{2}+\frac{1}{2})}{\Gamma(n+\frac{d}{2}+\frac{1}{2})\Gamma(n+1)}
 \notag \\
&=
\pi^{-\frac{1}{2}}
\frac{\Gamma (\frac{d}{2}+\frac{1}{2})}{\Gamma (\frac{d}{2})}
\sum\limits_{n=1}^{[\frac{t}{2}]}n^{-1}+\mathcal{O}(1)=
\pi^{-\frac{1}{2}}
\frac{\Gamma (\frac{d}{2}+\frac{1}{2})}{\Gamma (\frac{d}{2})}
\log t+\mathcal{O}(1) \notag \\
&=
\pi^{-\frac{1}{2}}
\frac{\Gamma (\frac{d}{2}+\frac{1}{2})}{d\Gamma (\frac{d}{2})}
\log N+\mathcal{O}(1).
\end{align}

Formulas (\ref{s_eq1}), (\ref{abs_valS_t})  and
(\ref{integrEst_s=d})
imply (\ref{theorems=d}) and (\ref{theorems=dLim}).
Theorem~\ref{theoremRisz} is proved.
\qed

\begin{acknowledgement}
The author is supported by the Austrian Science Fund FWF
  project F5503 (part of the Special Research Program (SFB) 
``Quasi-Monte Carlo Methods: Theory and Applications''
\end{acknowledgement}

%
%

\end{document}